\documentclass[11pt]{amsart}
\begin{document}
\title[Fredholm pairs and chains]{FURTHER RESULTS ON REGULAR FREDHOLM\\ PAIRS AND CHAINS}
\author[Enrico Boasso]{ENRICO BOASSO}
\maketitle
\footnotesize{The main objective of the present article is to characterize
regular Fredholm  pairs and\par
\indent chains in terms of Fredholm operators.\par
\par
}
\vskip.5truecm
\normalsize{
\centerline{\bf 1. INTRODUCTION}\par
\vskip.5cm
\indent In multiparameter operator theory there exist objects that in the context of Hilbert spaces
can be characterized using linear and bounded maps, to mention only some of the works related to
this subject, see for example [9, 10, 11, 12]. Since the adjoint plays a key role in this relationship between one and several
variable operator theory, it is not possible to directly
translate these results to the frame of Banach spaces. However in [7] R. Harte
and W. Y. Lee showed that generalized inverses could be used to reformulate the 
Hilbert space situation for regular Fredholm Banach space chains and complexes.\par

\indent On the other hand, Fredholm pairs and chains have been recently introduced and
their main properties have been studied, see [1, 2, 4, 5, 8] and the monograph [3]. These objects 
consist in generalizations of the notions of Fredholm operators and
Fredholm Banach space complexes respectively. In particular, 
in [4] Fredholm pairs and chains in Hilbert spaces were
characterized using the ideas mentioned in the previous paragraph.
Furthermore, in [5]
regular Fredholm pairs, that is Fredholm pairs whose operators
admit generalized inverses, were characterized and classified.\par

\indent The main objective of the present article is to characterize
regular Fredholm pairs and chains in Banach spaces using Fredholm operator
extending to these objects the characterizations
of [4] using the approach of R. Harte and W. Y. Lee in [7]. Naturally, as in [7], since 
the adjoint can not be considered any more, suitable generalized inverses
must be defined to prove the main results of this work. 
\par

\indent The article is organized as follows. In the next section the definitions
of the objects under consideration will be recalled. In section 3 two characterizations
of regular Fredholm pairs will be given. Finally, in section 4 the results
of section 3 will be applied to regular Fredholm chains to obtain two characterizations
of them.\par
\vskip.5cm
\centerline{\bf 2. PRELIMINARY DEFINITIONS }\par
\vskip.5cm
\indent From now on $X$ and $Y$ will denote two Banach spaces and $L(X,Y)$ will stand for the algebra of 
all linear and continuous operators defined on $X$ with values in $Y$. As usual, when $X=Y$, $L(X,X)$ is 
denoted by $L(X)$. For every $S\in L(X,Y)$, the null space of $S$ is denoted by
$N(S)=\{x\in X\colon \hbox{  }S(x)=0\}$, and the range of $S$ by $R(S)=
\{ y\in Y: \hbox{  }\exists
\hbox{ } x\in X \hbox{ such that }y=S(x)\}$. \par
\indent Next follow the definitions of Fredholm pairs and chains, 
see for instance [1, 3, 8].\par

\indent {\it Definition\rm} 2.1. Let $X$ and $Y$ be two Banach spaces and let $S\in L(X,Y)$
and $T\in L(Y,X)$ be such that the following dimensions are finite:
\begin{enumerate}
\renewcommand{\theenumi}{\roman{enumi}}
\item $a\colon =\dim N(S)/(N(S)\cap R(T))$, $b\colon = \dim R(T)/(N(S)\cap R(T))$,
\item $c\colon = \dim N(T)/(N(T)\cap R(S))$, $d\colon =\dim R(S)/(N(T)\cap R(S))$.
\end{enumerate}
A pair $(S,T)$ with the above properties is said to be a Fredholm pair.\par

\indent Let $P(X,Y)$ denote the set of all Fredholm pairs. If $(S,T)\in P(X,Y)$, then
the $index$ of $(S,T)$ is defined by the equality
\[
 \hbox{ind}  (S,T)\colon = a -b -c +d.
\]
\indent In particular, if $(S,T)\in P(X,Y)$ is such that $ST=0$
and $TS=0$, that is if $b=d=0$, then 
$(S,T)$ and $(T,S)$ are Fredholm chains in the sense of [6, section 10.6] and [7].\par

\indent {\it Definition\rm} 2.2. A Fredholm chain $(X_p,\delta_p)_{p\in Z}$ is 
a sequence of Banach spaces $X_p$ and bounded operators  $\delta_p\in L(X_p,X_{p-1})$ 
such that  there is a natural number $n$ with the property $X_p=0$, $p < 0$ and $p \ge n+1$,
$\delta_p=0$, $p \le 0$ and $p > n$, and
\[
N(\delta_p)/(N(\delta_p)\cap R(\delta_{p+1}))\hbox{ and }
R(\delta_{p+1})/(N(\delta_p)\cap R(\delta_{p+1}))
\] 
are finite dimensional subspaces of $X_p$, $p\in Z$.

\indent Given a Fredholm chain, it is possible
to associate an $index$  to it. In fact, if $(X_p,\delta_p)_{p\in Z}$ is such an object, then define 
\[
 \hbox{ind}  (X_p,\delta_p)_{p\in\mathit Z}=\sum_{p=0}^n(-1)^p d_p,
\]
where $d_p=\dim\hbox{ }
N(\delta_p)/(N(\delta_p)\cap R(\delta_{p+1}))-\dim R(\delta_{p+1})/(N(\delta_p)\cap R(\delta_{p+1}))$,
see [8].\par 

\indent Recall that in [8] the more general concept of 
semi-Fredholm chain was introduced.
However, since the main concern of this article consists in Fredholm objects,
only Fredholm chains will be considered. \par
\indent {\it Remark \rm } 2.3. There is a natural 
relationship between Fredholm pairs and chains. In fact, given as in Definition 2.2  a sequence of spaces and maps
$(X_p,\delta_p)_{p\in Z}$, 
consider the Banach spaces 
\[
X=\bigoplus_{p=2k}X_p,\hskip1cm Y=\bigoplus_{p=2k+1}X_p,
\]
and the Banach space operators  $S\in L(X,Y)$ and $T\in L(Y,X)$ defined  by 
\[
S=\bigoplus_{p=2k} \delta_p,\hskip1cm T=\bigoplus_{p=2k+1} \delta_p,
\]
where $X_p=0$, $p < 0$ and $p \ge n+1$, $\delta_p=0$, $p \le 0$ and $p > n$, and $n$ is 
a natural number.\par
\indent Now well, it is not difficult to prove that
dim $R(ST)$ and dim $R(TS)$ are finite dimensional if and only if
dim $R(\delta_{p+1})/(N(\delta_p)\cap R(\delta_{p+1}))$ are finite dimensional, p = 0, $\ldots$, n.
Furthermore, a straightforward calculation shows that necessary and sufficient  for $N(S)/(N(S)\cap R(T))$ and $N(T)/(N(T)\cap R(S))$ to be finite dimensional
is the fact that $N(\delta_p)/(N(\delta_p)\cap R(\delta_{p+1}))$ are finite dimensional, p = 0, $\ldots$, n.
Consequently, $(X_p,\delta_p)_{p\in Z}$ is a Fredholm chain if and only if $(S,T)$ is a Fredholm
pair. Finally, in this case,
\[
\hbox{ ind} (X_p,\delta_p)_{p\in Z}=\hbox{ ind} (S,T),
\]
see [4, Remark 2.4].\par

\indent On the other hand, recall that an operator $T\in L(X,Y)$ is called  $regular$, if there is $S\in L(Y,X)$
such that $T=TST$. The map $S$ is said to be a $generalized$ $inverse$ of $T$.
In addition, if  $T$ is also a generalized inverse of $S$, that is if $S=STS$,
then $S$ is said to be a $normalized$ $generalized$ $inverse$ of $T$. Note that if $T$ is regular,
then $T$ always has a normalized generalized inverse. In fact, if $S$ is a generalized
inverse of $T$, then $S'= STS$ is a normalized generalized inverse of $T$.
In the following definition, the notions of regular Fredholm pairs and chains will be recalled,
see [5].\par

\indent {\it Definition\rm} 2.4. Let $X$ and $Y$ be two Banach spaces and consider $S\in L(X,Y)$ and
$T\in L(Y,X)$ such that $(S,T)\in P(X,Y)$. The pair $(S,T)$ will be said to be a regular
Fredholm pair, if $S$ and $T$ are regular operators. Similarly, 
given a Fredholm chain $(X_p,\delta_p)_{p\in Z}$, then $(X_p,\delta_p)_{p\in Z}$ will be said to be a regular Fredholm chain,
if the Fredholm pair defined by $(X_p,\delta_p)_{p\in Z}$ is regular.\par

\indent Note that given $(S,T)\in P(X,Y)$, several statements equivalent to the fact that
$(S,T)$ is a regular Fredholm pair 
were considered in [5, Proposition 2.4].
Concerning regular Fredholm chains, see  [5, Remark 2.7].\par 
\vskip.5cm
\centerline{\bf 3. CHARACTERIZATIONS OF REGULAR FREDHOLM PAIRS}\par
\vskip.5cm
\indent First of all, a preliminary remark is presented.\par

\indent {\it Remark \rm }3.1. Let $X$ and $Y$ be two Banach spaces,
and let $S\in L(X,Y)$ and $T\in L(Y,X)$ be two operators such that
$R(ST)$ and $R(TS)$ are finite dimensional subspaces of $Y$ 
and $X$ respectively. Then, define the Banach spaces
$\mathcal{X}=X/R(TS)$ and $\mathcal{Y}=Y/R(ST)$,
and the linear and bounded maps $\tilde{S}\in L(\mathcal{X},\mathcal{Y})$
and $\tilde{T}\in L(\mathcal{Y},\mathcal{X})$
as the factorization of $S$ and $T$  through the respective invariant subspaces.
Clearly, $\tilde{S}\tilde{T}=0$ and $\tilde{T}\tilde{S}=0$,
that is  $(\tilde{S},\tilde{T})$ and  $(\tilde{T},\tilde{S})$ are chains in the sense of [7]. Moreover,
the regularity of $S$ and $T$ is equivalent to the one of  $\tilde{S}$ and $\tilde{T}$.\par

\indent THEOREM 3.2. \it  In the same conditions of Remark 3.1, the following statements are
equivalent.\par
\begin{enumerate}
\renewcommand{\theenumi}{\roman{enumi}}
\item $S$ and $T$ are regular operators.
\item $\tilde{S}$ and $\tilde{T}$ are regular operators.
\end{enumerate}
\rm
\indent {\it Proof. \rm }Since $S\in L(X,Y)$ is a regular operator, according to
 [6, Theorem 3.8.2], there is $M$ a linear subspace of $Y$ such that
$R(S)\oplus M=Y$. Consequently, $R(\tilde{S}) + \pi_Y (M)=\mathcal{Y}$,
where $\pi_Y\colon Y\to \mathcal{Y}$ is the quotient map.\par
\indent Next suppose that there is $\tilde{x}\in \mathcal{X}$ such that
$\tilde{S}(\tilde{x})=\tilde{m}\in \pi_Y (M)$. In particular,
there exist $x\in X$, $m\in M$, and $y\in Y$ such that $\pi_X(x)=\tilde{x}$,
 $\pi_Y(m)=\tilde{m}$, and $S(x)-m=ST(y)$, where $\pi_X\colon X\to \mathcal{X}$
is the quotient map.
However, in this case $S(x-T(y))=m$. Then, $m\in M\cap R(S)=0$, and
$\tilde{S}(\tilde{x})=\tilde{m}=0$. Therefore, $R(\tilde{S}) \oplus \pi_Y (M)=\mathcal{Y}$.\par
\indent On the other hand, according to [6, Theorem 3.8.2], there exists N,  a vector subspace
of $X$ such that $N(S)\oplus N=X$. However, since 
\[
(N(S)+R(T))/N(S)\cong R(T))/(N(S)\cap R(T))\cong R(ST),
\] 
there is a finite dimensional subspace $X_1$ such that
$N(S)\oplus X_1=N(S)+R(T)$. Then, a straightforward calculation
proves that there exists a closed vector subspace $R$ of $X$ such that
$(N(S)+R(T))\oplus R =X$.\par

\indent Now well, according to [1, Remark 2.1], $N(\tilde{S})=\pi_X(N(S)+R(T))$.
In particular, $N(\tilde{S})+\pi_X(R)=\mathcal{X}$. \par

\indent Next suppose that there are $n\in N(S)$, $y\in Y$, and $r\in R$
such that $\pi_X(n+T(y))=\pi_X(r)$. Then, there is $x_0\in X$ such that
$n+T(y)-r=TS(x_0)$. In particular, $n+T(y-S(x_0))=r$. Thus $r=0$,
which implies that $\pi_X (r)=\pi_X(n+T(x))=0$. Therefore,
$N(\tilde{S})\oplus\pi_X(R)=\mathcal{X}$. As a result,
according to what has been proved and to [6, Theorem 3.8.2], the operator $\tilde{S}$ is regular.\par

\indent Interchanging $S$ with $T$ and $\tilde{S}$ with $\tilde{T}$, it can be proved that $\tilde{T}$ is regular.\par

\indent In order to prove the converse implication, suppose that  
$\tilde{S}$ and $\tilde{T}$ are regular operators.\par

\indent Consider $\mathcal{V}$, a closed vector subspace of $\mathcal{Y}$, such that
$R(\tilde{S})\oplus \mathcal{V}=\mathcal{Y}$. Let $V_1=\pi_Y^{-1} (\mathcal{V})\cap R(ST)$.
Since $V_1$ has finite dimension, there is a vector subspace $W_1$ such
$V_1\oplus W_1=\pi_Y^{-1} (\mathcal{V})$. Moreover, since $V_1\subseteq R(ST)$ and
$\pi_Y$ is surjective, $\mathcal{V}=\pi_Y(W_1)$. Consequently,
$\pi_Y(R(S)+W_1)=R(\tilde{S})+ \mathcal{V}=\mathcal{Y}$,
which implies that $R(S)+W_1+R(ST)=Y$. However, since
$R(ST)\subseteq R(S)$, $R(S)+W_1=Y$.\par

\indent Next define $L=R(S)\cap W_1$. Then, $\pi_Y(L)\subseteq 
R(\tilde{S})\cap \mathcal{V}=0$. Thus, $L\subseteq R(ST)$. 
However, since $W_1\cap R(ST)=0$, $L=0$ and 
$R(S)\oplus W_1=X$.\par

\indent Similarly, suppose that $\mathcal{U}$ is a closed subspace of 
$\mathcal{X}$ such that $N(\tilde{S})\oplus \mathcal{U}=\mathcal{X}$.
Let $U_1=\pi_X^{-1} (\mathcal{U})\cap R(TS)$. Since $U_1$
is a finite dimensional subspace, there exists a subspace
$Z_1$ such that $U_1\oplus Z_1=\pi_X^{-1} (\mathcal{U})$.
Furthermore, since $\pi_X\colon X\to \mathcal{X}$
is a surjecive map, $\pi_X(Z_1)=\mathcal{U}$.\par

\indent Now well, according to [1, Remark 2.1],
$\pi_X(N(S)+R(T)+Z_1)=\pi_X(N(S) +R(T))+\pi_X (Z_1)
=N(\tilde{S}) +\mathcal{U}=\mathcal{X}$. Therefore,
since $R(TS)\subseteq R(T)$, $N(S)+R(T)+Z_1=X$.\par

\indent Let $P=(N(S)+R(T))\cap Z_1$. According to [1, Remark 2.1],  
$\pi_X(P)\subseteq N(\tilde{S})\cap \mathcal{U}=0$.
Thus, $P\subseteq R(TS)$. Consequently, since $Z_1\cap R(TS)=0$, $P=0$,
and $(N(S)+R(T))\oplus Z_1=X$. 
However, since $N(S)\oplus X_1= N(S)+R(T)$,
where $X_1$ is a finite dimensional vector subspace of $X$,
$N(S)\oplus (X_1+Z_1)=X$, which, according to what has
been proved and to [6, Theorem 3.8.2],
implies that $S$ is a regular operator.\par

\indent Interchanging $\tilde{S}$ with $\tilde{T}$ and $S$ with $T$,
it can be proved that $T$ is a regular operator.
\qed\par

\indent In order to prove the first characterization concerning regular Fredholm pairs, some preparation
is needed.\par

\indent {\it Remark \rm 3.3. Let $X$, $Y$, $S$, $T$,  $\mathcal{X}$, $ \mathcal{Y}$,
$\tilde{S}$, and $\tilde{T}$ be as in Remark 3.1. Recall that $\tilde{S}\tilde{T}=0$
and $\tilde{T}\tilde{S}=0$. Furthermore, according to Theorem 3.2,
necessary and sufficient for $S$ and $T$ to be regular bounded and linear maps
is the fact that $(\tilde{S},\tilde{T}, \tilde{S})$ is a regular chain in the sense of [7, p. 283-284].
In addition, according again to [7, p. 284],
there exist $\tilde{S}'\in L( \mathcal{Y}, \mathcal{X})$ and $\tilde{T}'\in L( \mathcal{X}, \mathcal{Y})$
such that $(\tilde{T}',\tilde{S}', \tilde{T}')$ is a regular chain in the sense of [7] with the property that
$\tilde{S}'$ (respectively $\tilde{T}'$) is a normalized generalized inverse of
$\tilde{S}$ (respectively $\tilde{T}$).\par
\indent Now well, since $R(ST)$ and $R(TS)$ are finite dimensional subspaces
of $Y$ and $X$ respectively, $\mathcal{X}$ and $\mathcal{Y}$  are isomorphic to
finite codimensional closed subspaces of $X$ and $Y$ respectively.
Consequently, using  Banach space isomorphisms, $\tilde{S}'$ and 
$\tilde{T}'$ can be thought of  operators defined on and to 
finite codimensional closed subspaces of $X$ and $Y$.
Donote by $S'\in L(Y,X)$ and $T'\in L(X,Y)$ the extension of $\tilde{S}'$ and $\tilde{T}'$ respectively,
such that $S'=0$ on
$R(TS)$ and $T'=0$ on $R(ST)$.
Note that  $S'$ and $T'$ are regular maps. \par

\indent Next follows the first characterization of regular Fredholm pairs.\par

\indent THEOREM 3.4. \it  In the same conditions of Remark 3.3, the following statements are equivalent.\par
\begin{enumerate}
\renewcommand{\theenumi}{\roman{enumi}}
\item $(S,T)$  is a regular Fredholm pair.
\item $S+T'\in L(X,Y)$ is a Fredholm operator.
\item $T+S'\in L(Y,X)$ is a Fredholm operator.
\end{enumerate}
\indent  Futhermore, in this case
\[
\hbox{ ind}  (S,T)=\hbox{  ind} (S+T') = -\hbox{ ind}  (T+S').
\]
\rm
\indent {\it Proof. \rm } According to [1, Remark 2.1] and to Theorem 3.2, $(S,T)$ is a regular Fredholm pair
if and only if $(\tilde{S},\tilde{T}, \tilde{S})$ is a Fredholm regular chain in the sense of [7].
Then, according to [7, Theorem 5], the first statement of the Theorem is
equivalent to the fact that $\tilde{S}+\tilde{T}'\in L(\mathcal{X},\mathcal{Y})$ is a Fredholm operator.\par

\indent Now well, since as in Remark 3.3 $\mathcal{X}$ and $\mathcal{Y}$ can be thought of
finite codimensional closed subspaces of $X$ and $Y$ respectively,
it is possible to define $S_1\in L(Y,X)$ as the extension of $\tilde{S}$
with the property $S_1=0$ on $R(TS)$. Consequently, since $R(ST)$ and $R(TS)$ are finite dimensional
vector subspaces, according to what has been proved, the first statement of the  Theorem
is equivalent to the fact that $S_1+T'\in L(Y,X)$ is a Fredholm operator.
However, since $S-S_1$ is a compact operator,  the first and the second statements of
the Theorem are equivalent.\par

\indent As regard the index, according to [1, Remark 2.1] and to [7, Theorem 5]
\[
\hbox{ ind} (S,T)-\dim R(TS) +\dim R(ST) = \hbox{ ind} (\tilde{S}, \tilde{T}) =\hbox{   ind} (\tilde{S}+\tilde{T}').
\]
\indent On the other hand, note that 
\[
\hbox{  ind} (S+T')= \hbox{  ind} (S_1+T') = \hbox{  ind} (\tilde{S}+\tilde{T}') +\dim R(TS) -\dim R(ST).
\]

\indent Therefore, $\hbox{ ind}(S,T) =\hbox{ ind} (S+T')$.\par
\indent A similar argument proves that the first and the third statements are equivalent as well as
the relationship between the indexes.
\qed\par

\indent The second characterization of regular Fredholm pairs is more general
in the sense that the generalized inverses can be chosen more freely, however, 
the formula regarding the index can not be considered. On the other hand,
to state this characterization, a Banach space operator is introduced.\par
 
\indent {\it Remark \rm 3.5  Let $X$, $Y$, $S$, $T$,  $\mathcal{X}$, $ \mathcal{Y}$,
$\tilde{S}$,  and $\tilde{T}$ be as in Remark 3.3. Since 
$\tilde{S}$,  and $\tilde{T}$ are regular operators, there exist
$\tilde{S}'\in L(\mathcal{Y},\mathcal{X})$ and  $\tilde{T}'\in L(\mathcal{X},\mathcal{Y})$  generalized inverses of $\tilde{S}$ and $\tilde{T}$ respectively.
In addition, proceeding as in Remark 3.3, denote by $S'\in L(Y,X)$ and $T'\in L(X,Y)$ any extensions of   
$\tilde{S}'$ and  $\tilde{T}'$ respectively.\par
\indent On the other hand, define the map $V\in L(X\oplus Y)$ as follows: $V\mid_ X \in L(X,Y)$, $V\mid_X=S+T'$, and $V\mid_Y\in L(Y,X)$,
$V\mid_Y=T+S'$.

\indent THEOREM 3.6. \it In the same conditions of Remark 3.5, the following statements are equivalent.\par
\begin{enumerate}
\renewcommand{\theenumi}{\roman{enumi}}
\item $(S,T)$  is a regular Fredholm pair.
\item $S'S+TT'$  and  $T'T + SS'$  are Fredholm operators. 
\item $V$  is a Fredholm operator.
\end{enumerate}
\rm
\indent {\it Proof. \rm } As in Theorem 3.4, according to [1, Remark 2.1] and to
Theorem 3.2, the first statement is satisfied if and only if 
$(\tilde{S},\tilde{T}, \tilde{S})$ is a Fredholm regular chain in the sense of [7],
which, according to [7, Theorem 3], is equivalent to  the fact that
$\tilde{S}'\tilde{S}+\tilde{T}\tilde{T}'\in L(\mathcal{X})$
and  $\tilde{T}'\tilde{T}+\tilde{S}\tilde{S}'\in L(\mathcal{Y})$
are Fredholm operators. Now well, since $R(TS)$ and $R(ST)$
are finite dimensional subspaces of $X$ and $Y$ respectively,
it is not difficult to prove that  what has been proved is equivalent
to the second statement of the Theorem.\par

\indent On the other hand, $V$ is a Fredholm operator if and only if $V^2$ also is Fredholm.
Now well, a straightforward calculation proves that there exists a finite range operator $F\in L(X\oplus Y)$ 
such that $V^2-F$ is a diagonal operator with entries $S'S+TT'\in L(X)$ and $SS'+TT'\in L(Y)$,
which clearly implies that the second and the third statements are equivalent.
\qed\par
\vskip.5cm
\newpage
\centerline{\bf 4. CHARACTERIZATIONS OF REGULAR FREDHOLM CHAINS }\par
\vskip.5cm
\indent In order to prove the main results of this section, some preliminary
facts will be considered.\par

\indent {\it Remark \rm 4.1 Let $(X_p,\delta_p)_{p\in Z}$ be a sequence of spaces and maps
such that $X_p=0$, $p < 0$ and $p \ge n+1$, and $\delta_p=0$, $p \le 0$ and $p > n$, where $n$ is 
a natural number. In addition, suppose that $R(\delta_p\delta_{p-1})$ is finite dimensional 
and $\delta_p\in L(X_p,X_{p-1})$ is a regular operator, $p\in Z$. Then, if
$X$, $Y$, $S\in L(X,Y)$, and $T\in L(Y,X)$ are as in Remark 2.3,
the properties of  $(X_p,\delta_p)_{p\in Z}$ are equivalent to the fact that
$R(ST)$ and $R(TS)$ are finite dimensional and $S$ and $T$ are
regular operators.\par
\indent Now well, if $\mathcal{X}$, $\mathcal{Y}$, $\tilde{S}$
and $\tilde{T}$ are as in Remark 3.3, note that
\[
\mathcal{X}=\bigoplus_{p=2k} \mathcal{X}_p, \hskip1.5truecm\mathcal{Y}=\bigoplus_{p=2k+1}\mathcal{X}_p,
\]
\[ 
\tilde{S}=\bigoplus_{p=2k}\tilde{\delta}_p, \hskip1.5truecm \tilde{T}=\bigoplus_{p=2k+1}\tilde{\delta}_p,
\]
where $\mathcal{X}_p=X_p/R(\delta_{p+1}\delta_{p+2})$ and $\tilde{\delta}_p\colon \mathcal{X}_p
\to \mathcal{X}_{p-1}$ is the quotient map induced by $\delta_p\colon X_p\to X_{p-1}$, $p\in  Z$.
Since $\tilde{S}\tilde{T}=0$ and $\tilde{T}\tilde{S}=0$, the sequence of spaces and maps
$(\mathcal{X}_p,\tilde{\delta}_p)_{p\in Z}$ is a finite complex of Banach spaces.
In addition, according to Theorem 3.2, $\tilde{S}$ and $\tilde{T}$ are regular operators,
which is equivalent to the fact that $\tilde{\delta}_p\colon \mathcal{X}_p
\to \mathcal{X}_{p-1}$ is regular, $p\in  Z$. A straightforward calculation
proves that there exists a sequence of linear and bounded maps $( \tilde{\delta}_p')_{p\in Z}$,
$\tilde{\delta}_p'\colon \mathcal{X}_{p-1}\to \mathcal{X}_p$ such that 
$\tilde{\delta}_p'$ is a normalized generalized inverse of $\tilde{\delta}_p$
and $(\mathcal{X}_p,\tilde{\delta}_p')_{p\in Z}$ is a finite complex of Banach spaces.
Consequently, if
\[
\tilde{S}'=\bigoplus_{p=2k}\tilde{\delta}_p',\hskip1.5truecm \tilde{T'}=\bigoplus_{p=2k+1}\tilde{\delta}_p',
\]
then $(\tilde{T}',\tilde{S}', \tilde{T}')$ is regular chain in the sense of [7]
such that $\tilde{S}'$ and $\tilde{T}'$ are normalized generalized inverses
of $\tilde{S}$ and $\tilde{T}$. Finally, as in Remark 3.3, consider $S'\in L(Y,X)$
and $T'\in L(X,Y)$, the extensions of $\tilde{S}'$ and $\tilde{T}'$ such that
$S'=0$ on $R(TS)$ and $T'=0$ on $R(ST)$,  and denote by $\delta_p'\colon X_{p-1} \to X_p$
the  extension of $\tilde{\delta}_p'$ such that
$\delta_p'=0$ on $R(\delta_p\delta_{p-1})$, $p\in Z$.

\indent THEOREM 4.2. \it In the same conditions of Remark 4.1, the following statements are equivalent.\par
\begin{enumerate}
\renewcommand{\theenumi}{\roman{enumi}}
\item $(X_p,\delta_p)_{p\in Z}$ is a regular Fredholm chain.
\item $\oplus_{p=2k} (\delta_p+\delta_{p+1}')$  is a Fredholm operator. 
\item $\oplus_{p=2k+1} (\delta_p+\delta_{p+1}')$  is a Fredholm operator.
\end{enumerate}

\indent Furthermore, in this case,
\[
\hbox{ ind} (X_p,\delta_p)_{p\in Z}= \hbox{ ind}\oplus_{p=2k} (\delta_p+\delta_{p+1}') =
-\hbox{ ind}\oplus_{p=2k+1} (\delta_p+\delta_{p+1}'). 
\]
\rm
\indent {\it Proof. \rm } Apply Remark 2.3, Theorem 3.4, and Remark 4.1.
\qed\par

\indent {\it Remark \rm } 4.3  In the same conditions of Remark 4.1, consider the same spaces and maps
but, as in Remark 3.5, $\tilde{S}'$ and $\tilde{T}'$ denote any generalized inverses of 
$\tilde{S}$ and $\tilde{T}$. In addition, $S'$ and $T'$ can be any extentions of  $\tilde{S}'$ and $\tilde{T}'$.
Similarly, $\delta_p'$ is any extension of $\tilde{\delta}_p'$, $p\in  Z$.\par

\indent THEOREM 4.4. \it  In the conditions of Remark 4.3, the following statements are equivalent.\par
 \begin{enumerate}
\renewcommand{\theenumi}{\roman{enumi}}
\item $(X_p,\delta_p)_{p\in Z}$  is a regular Fredholm chain.
\item $\delta_{p+1}\delta_{p+1}' +\delta_p'\delta_p\in L(X_p)$  is a Fredholm operator, $p\in Z$. 
\end{enumerate}
\rm
\indent {\it Proof. \rm } Apply Remark 2.3, Theorem 3.6, and Remarks 4.1 and 4.3.
\qed\par

\vskip.3truecm
\noindent Enrico Boasso
\par
\noindent E-mail address: enrico\_odisseo@yahoo.it
\end{document}